# Feasible Path SQP Algorithm for Simulation-based Optimization Surrogated with Differentiable Machine Learning Models

*Zixuan Zhang[1,2], Xiaowei Song[1,2], Yujiao Zeng[1], Jie Li[3,4], Yaling Nie[1], Min Zhu[1], Jianhua Chen[1], Linmin Wang[1,2], Xin Xiao,[1,4]\**

[1] Institute of Process Engineering, Chinese Academy of Science, Beijing, 100190, China

[2] School of Chemical Engineering, University of Chinese Academy of Science, Beijing 100049, China

[3] Centre for Process Integration, Department of Chemical Engineering, School of Engineering, The University of Manchester, Manchester M13 9PL, United Kingdom

[4] Carbon Neutral Intelligent Industry Innovation Center, China Institute for Innovation and Development Strategy, Beijing, 100044, China

**Abstract:**

With the development of artificial intelligence, simulation-based optimization problems, which present a significant challenge in the process systems engineering community, are increasingly being addressed with the surrogate-based framework. In this work, we propose a deterministic algorithm framework based on feasible path sequential quadratic programming for optimizing differentiable machine learning models embedded problems. The proposed framework effectively addresses two key challenges: (i) achieving the computation of first- and second-order derivatives of machine learning models' outputs with respect to inputs; and (ii) by introducing the feasible path method, the massive intermediate variables resulting from the algebraic formulation of machine learning models eliminated. Surrogate models for six test functions and two process simulations were established and optimized. All six test functions were successfully optimized to the global optima, demonstrating the framework's effectiveness. The optimization time for all cases did not exceed 2s, highlighting the efficiency of the algorithm.

---

\* Corresponding author: Prof. Xin Xiao, E-mail: xxiao@ipe.ac.cn.

**Keywords:** simulation-based optimization; surrogate model; machine learning; feasible path algorithm; sequential quadratic programming

## 1. Introduction

Optimizing the unit or factory with strict simulation model is an essential task for the development of chemical industry[1]. According to the characteristics of the simulator, the optimization problems with strict simulation models can be divided into three categories: sequential module (SM)-based problem, equation-oriented (EO)-based problem and surrogate-based problem. SM method and EO method are two main mechanism-driven simulation methods. However, the SM method is arduous to obtain the accurate derivatives, and the EO method will create enormous number of variables and constraints[2,3], both making optimizing particularly difficult. Another modeling approach uses a data-driven model to construct a mapping between input and output through a large amount of data generated by the mechanism model or the real-world factories.

To address problems that rely on SM simulators, derivative-free optimization (DFO) methods are commonly employed[4,5]. These techniques can be broadly divided into two approaches: search-based such as particle swarm optimization (PSO) algorithm or DIRECT (DIviding RECTangles) algorithm[6], and model-based like trust-region methods[7,8] and kriging methods[9], who establish an approximate algebraic model using a tiny amount of simulated data in a local area. Javaloyes-Antón et al.[10] and Lee et al.[11] have successfully applied heuristic algorithm to solve optimization problems with embedded SM simulators. Na et al.[5] use deterministic global search algorithm DIRECT to optimize liquefied natural gas (LNG) process with hidden constraints. However, search-based approaches require the repeated computation of expensive objective function to obtain a certain level of global search capability, and this can be extremely time-consuming. The model-based approach leverages a minimal set of local sampling points to construct an algebraic model, whose derivatives can facilitate the update of search point. The response surface method, represented by Kriging models[12,13], is one of the earliest to be used in model-based algorithms. Cozad et al.[14] then proposed a framework named automated learning of algebraic models for optimization (ALAMO), who fit the simulation model with several basis functions and modify the model in optimizing iterations. Ma et al.[15] proposed an Branch-and-Model (BAM) method to augment the global optimality of ALAMO. Eason and Biegler[7,8] present a promising method combined with trust-region filter sequential quadratic programming (SQP) algorithm and the idea of derivative-free optimization to solve glass/black-box problems. While model-based approach demands less computational effort compared to search-based techniques, it may also lead to local optima, constraint violations, and non-convergence.

EO-based problems, taking the closed nonlinear equation system that the simulator needs to solve as a set of constraints, are usually formulated as large nonlinear programming (NLP). EO model allows the use of powerful decomposition methods such as Lagrange decomposition and Schur-complement to break down large-scale NLP problems into smaller and more easily solvable subproblems[16]. This strategy not only simplifies the computational complexity but also enhances the efficiency of the optimization process by leveraging the inherent structure of the problem. Kang et al.[17] use preconditioned conjugate gradient (PCG) method to solve the Schur-complement system and apply this approach in optimization of large nonlinear systems. Decomposition methods, such as orthogonal decomposition[18–20], can also be employed to establish reduced-order models, representing another widely adopted strategy for EO-based problems. By constructing a lower-dimensional representation of the original system, reduced-order models are capable of capturing the essential features and their relationship while significantly reducing the computational burden[21–23]. Different from SM model, EO models are capable of providing Jacobian matrix for optimization, which has greatly catalyzed the application of feasible path algorithm in solving EO-based problems. Ishii and Otto[24] first introduced a two-tier SQP-based decomposition algorithm, where the inner tier is tasked with calculating the numerical solution and the Jacobian matrix of the EO model, and the outer tier is responsible for updating the search point and dealing with the constraints. This structured approach is collectively termed as the feasible path method[25–27] in process system engineering (PSE) community. Ma et al.[27] substituted the steady-state simulation model[28,29] with pseudo-transient continuation model[30–32], which guarantees the convergence of the EO simulator. Subsequently, they improved SQP and sequential least squares programming (SLSQP) to bolster the algorithm's robustness[26,33]. Liu et al.[34,35] employed this method in branch and bound algorithm to solve EO-based MINLP. Zhao et al.[2] and Zhang et al.[36] use symbolic computation rather than numerical computation to solve nonlinear system, attempting to obtain the analytical solution to accelerate the simulation, which can also be considered as a specific form of feasible path method. Despite advancements already made, the methods mentioned above remain reliant on EO simulators or large-scale NLP algorithms, both consuming substantial computation costs in each iteration.

With the development of artificial intelligence (AI) technology, surrogate models, especially machine learning (ML) models, provide a rapid computing speed and strong regression ability, which can effectively compensate for the difficulties in optimization involving rigorous models, and now has been widely adopted by PSE community. Distinguished from model-based algorithm of DFO, surrogate-based problem establishes the approximate models with a larger dataset before optimizing, rather than sampling, learning and modifying them in each iteration. The essence of surrogate model is to construct a mapping between the decision variables and the dependent variables that are relevant to the objective function and constraints. A novel and powerful regression tool is symbolic regression who is composed of a tree structure and multiple operators[37], and allowed to be expressed in the form of composite functions. Both genetic programming[38–40] (GP) and MINLP[41] can be employed to learn its representation. Furthermore, Kolmogorov–Arnold networks[42] (KAN) can also be regarded as an advanced rendition of symbolic regression, pushing the boundaries of what is achievable in terms of modeling complexity and interpretability. Although ML models are usually regarded as black-box models, numerous scholars apply their explicit algebraic expressions to deterministic solvers[43–49] with approaching the resolution of surrogate-based problems by employing the explicit formulations of neural networks as nonlinear constraints[47–49]. Concurrently, another cohort of researchers has explored the utilization of decision trees[50,51] or neural networks with rectified linear unit (ReLU) activations[43,52–55] who can be reformulated as mixed integer constraints. This strategy is aimed at tackling the surrogate-based problem with mixed integer linear programming (MILP)[51,54] or transforming nonconvex NLP into MILP which can be globally solved[50,55]. However, this reformulated optimization problem based on explicit format of ML models may be time-consuming due to a large number of exponential terms or binary variables, which undermines the original intention of taking advantage of the fast computing speed of ML models.

In this paper, we intended to develop a deterministic algorithm for surrogate-based problems, that not only preserves the fast computation speed of surrogate models, but also can quickly obtain the accurate first- and second-order derivative. To accomplish our ambition, the algorithm is divided into internal and external layers. The SQP-driven feasible path algorithm which is widely used in EO-based problems[26,33], serving as the external layer, takes charge of updating the direction and step size of the decision variables, who is also taken as the inputs of the surrogate model. In internal layer, the forward and backward of the ML model will be executed, so that the value of the state variables, who is also taken as the outputs of the surrogate model, and derivative information will be obtained, which can guide the update of decision variables in the external layer.

The main contributions of this work are as follows:

1) We present a generalized deterministic algorithmic framework based on SQP-driven feasible path method that permits the ML models embedded, offering both theoretical guarantees and rapid computational speed;

2) We derived the explicit algebraic expression of the ML models and both the first- and second-order derivative of the output with respect to the input, ensuring that quadratic programming (QP) subproblem can be generated in proposed algorithm workflows;

3) We conduct six test functions and two case studies to demonstrate the effectiveness and efficiency of the proposed framework.

The structure of the article is as follows. The overall surrogate-based optimization problem is defined in Section 2. The explicit algebraic expression and derivatives of common ML models is derived in Section 3. Section 4 offers the generalized framework based on SQP-driven feasible path algorithm. The computational experiments of cases are conducted, analysed and discussed in Section 5. Finally, Section 6 provides conclusion and prospect.

## 2. Problem description

In this work, we attempt to solve a surrogate-based optimization problem (P) with the following structure:

$$\begin{aligned} \min_{\mathbf{x}\in\mathbb{R}^n, \mathbf{y}\in\mathbb{R}^m} \quad & f(\mathbf{x},\mathbf{y}) \\ \text{s.t.} \quad & \mathbf{y}-\mathbf{s}(\mathbf{x})=0 \\ & \mathbf{h}(\mathbf{x},\mathbf{y})=0 \\ & \mathbf{g}(\mathbf{x},\mathbf{y})\leq 0 \end{aligned} \quad (P)$$

Where, the function $f: \mathbb{R}^n \times \mathbb{R}^m \to \mathbb{R}$ is the objective function, $\mathbf{s}: \mathbb{R}^n \to \mathbb{R}^m$ is the explicit algebraic expression of the surrogate model as is presented in Section 3, $\mathbf{h}: \mathbb{R}^n \times \mathbb{R}^m \to \mathbb{R}^{p_E}$ is the equality constraints and $\mathbf{g}: \mathbb{R}^n \times \mathbb{R}^m \to \mathbb{R}^{p_I}$ is the inequality constraints, $p_E$ and $p_I$ represent the dimensions of equations and inequalities; vector variables $\mathbf{x}$ and $\mathbf{y}$ represent decision variables and state variables, respectively, and are also the inputs and outputs of the surrogate model.

According to the feasible path method, we consider $\mathbf{y}-\mathbf{s}(\mathbf{x})=0$ as a complete system of equations, and $\mathbf{y}$ can be expressed by $\mathbf{x}$. Therefore, we can replace the variables $\mathbf{y}$ with $\mathbf{x}$, and reformulate problem (RP) as follows:

$$\begin{aligned} \min_{\mathbf{x}\in\mathbb{R}^n} \quad & f(\mathbf{x},\mathbf{s}(\mathbf{x})) \\ \text{s.t.} \quad & \mathbf{h}(\mathbf{x},\mathbf{s}(\mathbf{x}))=0 \\ & \mathbf{g}(\mathbf{x},\mathbf{s}(\mathbf{x}))\leq 0 \end{aligned} \quad (RP)$$

With this strategy, we can calculate $\mathbf{s}(\mathbf{x})$, $\nabla_\mathbf{x}\mathbf{s}(\mathbf{x})$ and $\nabla_\mathbf{x}^2\mathbf{s}(\mathbf{x})$ in the internal layer, which is used for updating objective function and variables $\mathbf{x}$, and only $\mathbf{x}$ need to be updated in the external layer.

## 3. Differentiable machine learning models

Surrogate models are extensively utilized in engineering design and optimization[56–60].

Simple surrogate models, such as polynomial models[61,62] and piecewise linear models, can be readily integrated into the optimization process. However, due to their straightforward structure, these models often exhibit low accuracy and poor generalization capabilities. This implies that when addressing complex or high-dimensional problems, these models may fail to accurately represent the system's behavior, leading to suboptimal solutions.

Recently, ML models have gained widespread application in process modeling owing to their strong fitting and extrapolation abilities. For instance, models such as neural networks[43,52,53], support vector machines[55], and decision tree[50,51] can capture nonlinear relationships within data, providing more precise predictions. However, the intricate structures of these models make them challenging to directly incorporate into derivative-based optimization algorithms. A common approach to integrating machine learning models with derivative-based algorithms is to embed the algebraic form of the ML model as constraints within the optimization model[43,52,63–66]. However, this practice often introduces a large number of non-convex terms[53,54], making the optimization problem significantly more difficult to solve.

Consequently, when employing ML models for optimization, heuristic algorithms or deterministic derivative-free algorithms are commonly used. These algorithms do not depend on gradient information but rather optimize through direct comparison of function values. However, these algorithms generally converge more slowly and incur higher computational costs when dealing with large-scale or high-dimensional problems.

Derivative-based algorithms possess numerous significant advantages, including faster convergence speed and the ability to handle complex constraints. These features render derivative-based algorithms highly popular in various practical applications. However, to integrate machine learning models into derivative-based optimization algorithms, the issue of derivative computation must be addressed.

In this section, we will introduce several commonly used ML models and present the algebraic expressions along with their derivatives. Rather than treating the algebraic expression of the ML model as a constraint and computing its Jacobian matrix, we will preserve the integrity of the ML model by using it as a composite function for differentiation. This approach embodies the essence of the feasible path method. Given that SQP is essentially equivalent to the Newton method[67], it is crucial to derive both the first-order derivative $\nabla_{\mathbf{x}} s(\mathbf{x})$ and the second-order derivative $\nabla_{\mathbf{x}}^2 s(\mathbf{x})$, which are used for formulating QP subproblems. To ensure mathematical rigor, the continuity and differentiability requirements of the ML models must be carefully considered. However, requiring that $s(\mathbf{x})$ belongs to $\mathcal{C}^2$, meaning it needs to be both continuous and twice differentiable, may be overly restrictive for surrogate models. Consequently, it is sufficient to demand the ML model to meet any of the following criteria within its domain:

1) $s(\mathbf{x})$ belongs to $\mathcal{C}^2$ and $\nabla_\mathbf{x}^2 s(\mathbf{x}) \neq \mathbf{0}$;

2) $s(\mathbf{x})$ belongs to $\mathcal{C}^0$, the basis functions $p(\cdot)$ of the ML model belong to $\mathcal{C}^2$, and $\nabla^2 p \neq \mathbf{0}$.

## 3.1 Decision Tree (DT)

Various forms of decision trees have been widely used as surrogate models in optimization[45,46,51,63]. In the early stage, the decision boundaries of decision trees were univariate and parallel to the coordinate axes, using the mean as the predictive model. Subsequently, decision trees with linear predictive models[46,51] were developed. In recent studies, decision trees can employ linear hyperplanes as decision boundaries and use arbitrary functions as predictive models to fit sample data[68].

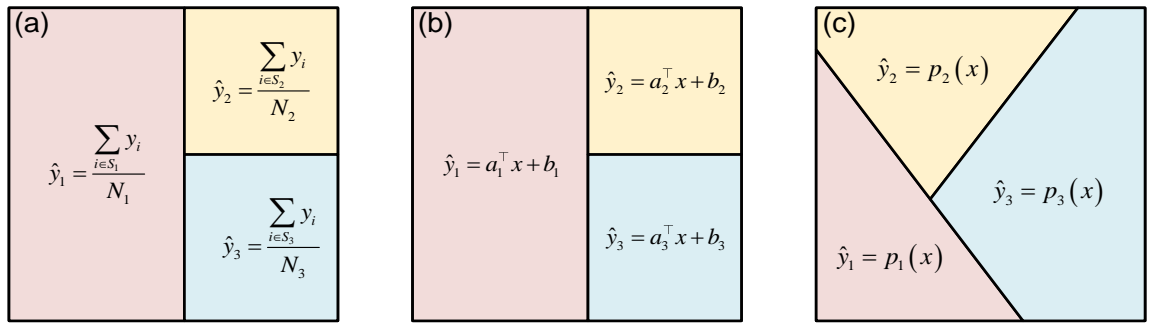

FIGURE 1 Decision tree with (a) single variable decision boundary and mean prediction model, (b) single variable decision boundary and linear prediction model as well as (c) linear hyperplane decision boundary and arbitrary prediction model.

In this work, we do not make any requirement to the form of the decision boundary, but require the predictive model in each leaf is quadratic differentiable and the second-order derivative is not equal to zero. Considering the single variable decision boundary can be regarded as a special form of linear decision boundary, the generalized algebraic formulation can be state as Eqs.(3) and (4).

$$\hat{y}_{DT} = p_i(x) \quad x \in \mathcal{P}_i, i \in \mathcal{L} \qquad (3)$$

$$\mathcal{P}_i = \left( \bigcap_{j \in L(L_i)} \{x \in \mathbb{R}^n : a_j^\top x \leq b_j\} \right) \cap \left( \bigcap_{j \in R(L_i)} \{x \in \mathbb{R}^n : a_j^\top x > b_j\} \right) \qquad (4)$$

Where, $\mathcal{L}$ is the set of leaves; $\hat{y}_{DT}$ is the predictive value of decision tree, $x$ is the input feature; $L_i$ represents the leaf node $i$, $j$ represents the branch node; $p_i(\cdot)$ is the the predictive model at $L_i$. For a binary tree, the path from the root to the leaves exists and is unique. For branch node $j$ in the path of leaf $L_i$, as is shown in FIGURE 2 (a), if $x$ satisfies $a_j^\top x \leq b_j$, then the node $j$ belongs to $L(L_i)$, otherwise belongs to $R(L_i)$. $\mathcal{P}_i$ represents the polyhedron as is stated in Eq.(4).

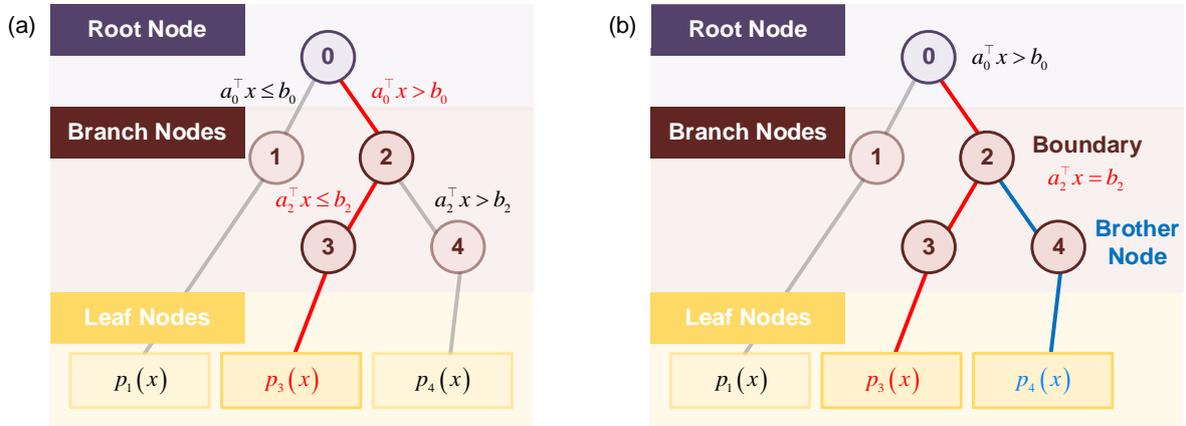

FIGURE 2 (a) The path from the root to the leaves and (b) the boundary prediction functions' recognition.

For those points that do not exactly satisfy the equality boundary, the first- and second-order derivatives are as follows:

$$\nabla_x \hat{y}_{DT} = \nabla_x p_i(x) \quad x \in \left( \bigcap_{j \in L(L_i)} \{a_j^\top x < b_j\} \right) \cap \left( \bigcap_{j \in R(L_i)} \{a_j^\top x > b_j\} \right) \qquad (5)$$

$$\nabla_x^2 \hat{y}_{DT} = \nabla_x^2 p_i(x) \quad x \in \left( \bigcap_{j \in L(L_i)} \{a_j^\top x < b_j\} \right) \cap \left( \bigcap_{j \in R(L_i)} \{a_j^\top x > b_j\} \right) \quad (6)$$

For those points that exactly satisfy the equality boundary, they are non-differentiable. Therefore, the minimum norm subgradient is used in place of the first-order derivative:

$$\partial \hat{y}_{DT} = \min_{\omega_i} \left\| \sum_i \omega_i \nabla_x p_i(x) \right\|_2 \quad x \in \bigcup_{j \in L(L_i) \cup R(L_i)} \{a_j^\top x = b_j\} \quad (7)$$

Then, employ the finite difference method to estimate the gradients of subgradient and construct an approximate Hessian matrix. The algorithm for boundary nodes' identification and derivatives' calculation of decision tree is shown in ALGORITHM 1.

ALGORITHM 1 Boundary nodes' identification and derivatives' calculation of decision tree.

| | |
|---|---|
| Step1: | Identify the path from root to leaf that $x$ belongs to. And determine whether $x$ satisfies the equality boundary of the branch node. |
| Step2: | If no equation is satisfied, proceed to Step3. Otherwise, proceed to Step4. |
| Step3: | Use Eqs.(5) and (6) to calculate derivatives and proceed to terminal. |
| Step4: | Take the brother node of boundary node as the root of the subtree, and perform Step1 on the subtree. If no equation is satisfied, proceed to Step5. Otherwise, repeat Step4. |
| Step5: | Use Eq.(7) and finite difference to calculate derivatives and proceed to terminal. |

## 3.2 Support Vector Regression Machine (SVR)

SVR is the application of Support Vector Machine (SVM) in regression problems. Similar to finding an optimal hyperplane that maximizes the margin between different classes in SVM, SVR seeks to identify a function that best fits the data while maintaining a certain tolerance for error, achieved through the introduction of an ε-insensitive loss function[69]. This model is especially effective for handling high-dimensional data and can construct non-linear decision boundaries in the feature space when kernel methods are utilized.

The algebraic form of kernel SVR is shown in Eq.(8):

$$\hat{y}_{SVR} = \sum_i (\alpha_i - \alpha_i^*)\mathcal{K}(x, x_i) + \beta \tag{8}$$

Where $\hat{y}_{SVR}$ is the predictive value of SVR, $\alpha_i$ and $\alpha_i^*$ is the Lagrange multiplier, $\beta$ is the bias term, $x_i$ is the sample $i$ in the training set, $x$ is the input feature and $\mathcal{K}(\cdot)$ is the kernel function.

And the derivatives are as shown in Eqs.(9) and (10).

$$\nabla_x \hat{y}_{SVR} = \sum_i (\alpha_i - \alpha_i^*)\nabla_x \mathcal{K}(x, x_i) \tag{9}$$

$$\nabla_x^2 \hat{y}_{SVR} = \sum_i (\alpha_i - \alpha_i^*)\nabla_x^2 \mathcal{K}(x, x_i) \tag{10}$$

## 3.3  Neural Networks (NN)

Neural networks serve as the cornerstone of deep learning technology, and their applications within chemical engineering are expanding rapidly. This growth is primarily attributed to neural networks' robust nonlinear modeling capabilities and their proficiency in processing intricate data.

By leveraging these strengths, neural networks enable advanced process optimization, predictive maintenance, and innovative product development in the chemical industry, thereby transforming traditional practices with more efficient and intelligent solutions.

The most basic NN model is the multi-layer perceptron (MLP). The algebraic form of MLP with $L$ hidden layers can be expressed by the following composite function:

$$\hat{y}_{MLP} = f^{(L+1)} \circ \phi^{(L)} \circ f^{(L)} \circ \cdots \circ \phi^{(1)} \circ f^{(1)}(x) \tag{11}$$

Where, $\hat{y}_{MLP}$ is the predictive value of MLP, $f^{(n)}(z) = W^{(n)}z + b^{(n)}$, $\phi^{(n)}(\bullet)$ is the activation function of the n-th layer. $W^{(n)}$ and $b^{(n)}$ represent the weights and biases from layer $n-1$ to layer $n$, respectively.

According to the chain differentiation rule, we have:

$$\frac{\partial \hat{y}_{MLP}}{\partial x} = \frac{\partial f^{(L+1)}}{\partial \phi^{(L)}} \cdot \frac{\partial \phi^{(L)}}{\partial f^{(L)}} \cdots \frac{\partial \phi^{(1)}}{\partial f^{(1)}} \cdot \frac{\partial f^{(1)}}{\partial x} \tag{12}$$

Specifically, for every layer we have:

$$\frac{\partial f^{(n+1)}}{\partial \phi^{(n)}} = W^{(n+1)} \tag{13}$$

$$\frac{\partial \phi^{(n)}}{\partial f^{(n)}} = diag\left(\phi'^{(n)}\left(f^{(n)} \circ \phi^{(n-1)} \circ f^{(n-1)} \circ \cdots \circ \phi^{(1)} \circ f^{(1)}(x)\right)\right) \tag{14}$$

Therefore, the first-order derivative can be expressed as follows:

$$\frac{\partial \hat{y}_{MLP}}{\partial x} = W^{(L+1)} \prod_{n=L}^{1} diag\left(\phi'^{(n)}\left(f^{(n)} \circ \phi^{(n-1)} \circ f^{(n-1)} \circ \cdots \circ \phi^{(1)} \circ f^{(1)}(x)\right)\right) W^{(n)} \tag{15}$$

To simplify the expression of second-order derivative of MLP, we assume that the predicted output value $\hat{y}_{MLP}$ is a scalar, and the conclusions obtained can be easily extended to high-dimensional output neural networks.

From the output layer to the hidden layer $L$, we have:

$$\frac{\partial^2 \hat{y}_{MLP}}{\partial \mathbf{x} \partial \mathbf{x}^\top} = W^{(L+1)} \cdot \frac{\partial^2 \phi^{(L)}}{\partial \mathbf{x} \partial \mathbf{x}^\top} \tag{16}$$

For hidden layer $l$ we have:

$$\frac{\partial^2 \phi^{(l)}}{\partial \mathbf{x} \partial \mathbf{x}^\top} = \phi''^{(l)}\left(f^{(l)}\right) \cdot \left(\frac{\partial f^{(l)}}{\partial \mathbf{x}}\right) \cdot \left(\frac{\partial f^{(l)}}{\partial \mathbf{x}}\right)^\top + \phi'^{(l)}\left(f^{(l)}\right) \cdot \frac{\partial^2 f^{(l)}}{\partial \mathbf{x} \partial \mathbf{x}^\top} \quad (17)$$

Where,

$$\frac{\partial f^{(l)}}{\partial \mathbf{x}} = \prod_{m=1}^{l-1}\left(diag\left(\phi'^{(m)}\left(f^{(m)}\right)\right)W^{(m)}\right) \quad (18)$$

$$\frac{\partial^2 f^{(l)}}{\partial \mathbf{x} \partial \mathbf{x}^\top} = \sum_{m=1}^{l-1}\left[\left(\prod_{r=m+1}^{l-1}\left(diag\left(\phi'^{(r)}\left(f^{(r)}\right)\right)W^{(r)}\right)\right) \cdot \frac{\partial^2 \phi^{(m)}}{\partial \mathbf{x} \partial \mathbf{x}^\top} \cdot \left(\prod_{k=1}^{m-1}\left(diag\left(\phi'^{(k)}\left(f^{(k)}\right)\right)W^{(k)}\right)\right)\right]$$

(19)

By applying Eqs.(16)-(19) and performing layer-by-layer recursion, we can obtain the second-order derivative of the MLP. Obviously, the algebraic symbolic differentiation approach yields expressions that are overly complex and challenging to implement. Therefore, we provide an alternative method: reverse-mode automatic differentiation[70,71] (AD). Different from central differencing and other numerical differentiation techniques, automatic differentiation provides exact derivatives while significantly reducing computational complexity.

Assuming we have a neural network as shown in FIGURE 3 (a), and the forward process can be transformed into a computational graph as is shown in FIGURE 3 (b). Here, $x_m$ is the input feature of NN, $y$ is the prediction value. $v_i$ and $v_k$ represent the linear transformation in NN, for example $v_1 = w_1 x_1$, where $w$ represent the weights. $p_j$ represents the activation function, for instance $p_1 = \phi(v_1 + v_4 + b)$, where $b$ represents the bias. By transforming NN into computational graph, the forward trace can be distinctly presented and the backward trace is used for calculating first-order derivative. Similarly, this process can also be transformed into a computational graph used for calculating second-order derivatives. In this work, we use PyTorch[72] to easily achieve automatic differentiation without the tedious derivation for algebraic expressions.

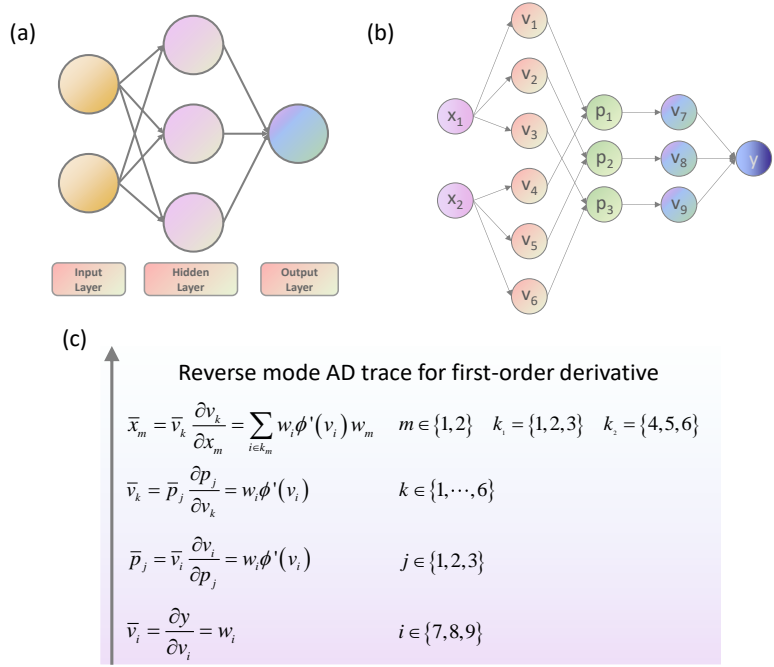

FIGURE 3 An illustrative example of NN with its (a) structure, (b) computational graph and (c) reverse mode AD trace for first-order derivative.

## 3.4 Ensemble Learning

For some special tasks, a single learner may not perform well. Therefore, the multiple base learners need to be integrated. The weighting method is the most commonly used ensemble approach in regression tasks:

$$\hat{y}_E = \sum_i \omega_i \hat{y}_i \qquad (20)$$

Where, $\hat{y}_E$ is the prediction value of ensemble learner, $\hat{y}_i$ is the prediction value of base learners and $\omega_i$ is the weights. The derivatives of integrated learner are also the weighted derivatives of the base learner.

The methods outlined above enable the provision of necessary derivative information while preserving the powerful fitting and rapid inference capabilities of machine learning models, thereby facilitating efficient optimization.

## 4. SQP-Driven feasible path algorithm

In this section, we will introduce SQP algorithm for problem (RP). Once the derivative information obtained, the SQP algorithm can be easily applied to simulation-based optimization problems surrogated with ML models. The Lagrangian function of (RP) is shown as Eq.(21):

$$L(\mathbf{x},\mathbf{s}(\mathbf{x}),\boldsymbol{\lambda},\boldsymbol{\mu}) = f(\mathbf{x},\mathbf{s}(\mathbf{x})) + \boldsymbol{\lambda}^\top \mathbf{h}(\mathbf{x},\mathbf{s}(\mathbf{x})) - \boldsymbol{\mu}^\top \mathbf{g}(\mathbf{x},\mathbf{s}(\mathbf{x})). \tag{21}$$

Where, $\boldsymbol{\lambda}$ and $\boldsymbol{\mu}$ are the Lagrange multipliers for equality and inequality constraints, respectively.

In this work, we take $\mathbf{H}^k := \nabla_{\mathbf{x}}^2 L(\mathbf{x}^k, \mathbf{s}(\mathbf{x}^k), \boldsymbol{\lambda}^k, \boldsymbol{\mu}^k)$ as the Hessian matrix of the Lagrangian function at $\mathbf{x}^k$, and $\mathbf{B}^k := \mathbf{H}^k + \mathbf{E}^k$ as the modified Hessian matrix. $\mathbf{E}^k$ denotes the minimum modification under the Frobenius norm. In order to obtain $\mathbf{E}^k$, first perform spectral decomposition on $\mathbf{H}^k$, so that $\mathbf{H}^k = \mathbf{Q}\boldsymbol{\Lambda}\mathbf{Q}^\top$, where $\mathbf{Q}$ is a matrix composed of eigenvectors of $\mathbf{H}^k$ and $\boldsymbol{\Lambda}$ is a diagonal matrix composed of eigenvalues of $\mathbf{H}^k$. In this way, we take $\delta$ as a constant approaching 0, and by modifying the eigenvalues of $\mathbf{H}^k$, we adjust all eigenvalues who are smaller than $\delta$ to $\delta$, thus ensuring the positive definiteness of $\mathbf{B}^k$. Therefore, we can define $\mathbf{E}^k := \mathbf{Q}\mathrm{diag}(\tau_i)\mathbf{Q}^\top$, where

$$\tau_i = \begin{cases} 0 & \sigma_i > \delta \\ \delta - \sigma_i & \sigma_i < \delta \end{cases},$$

and $\sigma_i$ represents the eigenvalues of $\mathbf{H}^k$.

## 4.1 SQP Algorithm

In this work, we use SQP to solve an optimization problem embedded with a differentiable ML model, who expands the original problem into a convex quadratic programming problem through Taylor expansion at each iteration point and conduct line search to merit function along the direction obtained from problem (QP). SQP is a type of Newton's method that can effectively handle both small and large-scale nonlinear constrained optimization problems. It can be integrated into line search or trust region frameworks. In this work, we use the line search SQP. Perform second-order Taylor expansion on the Lagrangian function of (RP) at $\mathbf{x}^k$ and first-order Taylor expansion on the constraints to obtain the following quadratic programming problem (QP):

$$\begin{aligned} \min_{\mathbf{x}\in\mathbb{R}^n} \quad & \frac{1}{2}\mathbf{d}^{k\top}\mathbf{B}^k\mathbf{d}^k + \nabla f\left(\mathbf{x}^k,\mathbf{s}(\mathbf{x}^k)\right)^\top \mathbf{d}^k \\ \text{s.t.} \quad & \nabla \mathbf{h}\left(\mathbf{x}^k,\mathbf{s}(\mathbf{x}^k)\right)^\top \mathbf{d}^k + \mathbf{h}\left(\mathbf{x}^k,\mathbf{s}(\mathbf{x}^k)\right) = 0 \\ & \nabla \mathbf{g}\left(\mathbf{x}^k,\mathbf{s}(\mathbf{x}^k)\right)^\top \mathbf{d}^k + \mathbf{g}\left(\mathbf{x}^k,\mathbf{s}(\mathbf{x}^k)\right) \leq 0 \end{aligned} \quad (\text{QP})$$

Where, $\mathbf{d}^k = (\mathbf{x} - \mathbf{x}^k)$ is the search direction of (RP), and $\nabla f$, $\nabla \mathbf{h}$ and $\nabla \mathbf{g}$ are the gradients of $f$, $\mathbf{h}$ and $\mathbf{g}$ respectively. The derivative with respect to $x_i$ is

$$\nabla_{x_i} F(\mathbf{x}, \mathbf{s}(\mathbf{x})) = \frac{\partial F}{\partial x_i} + \sum_{j=1}^{m} \frac{\partial F}{\partial s_j(\mathbf{x})} \frac{\partial s_j(\mathbf{x})}{\partial x_i},$$

where $F$ can be any one of $f$, $h$ and $g$.

After solving the problem (QP), we obtain the search direction $\mathbf{d}^k$, multipliers $\boldsymbol{\lambda}_{qp}^k$ and $\boldsymbol{\mu}_{qp}^k$ of problem (QP) at $\mathbf{x}^k$. Then, line search will be executed on the merit function and the step size $\alpha$ will be determined through backtracking method. In next iteration, $\mathbf{x}^{k+1} = \mathbf{x}^k + \alpha \mathbf{d}^k$ and $\mathbf{H}^{k+1} := \nabla_{\mathbf{x}}^2 L(\mathbf{x}^{k+1}, \mathbf{s}(\mathbf{x}^{k+1}), \boldsymbol{\lambda}_{qp}^k, \boldsymbol{\mu}_{qp}^k)$, according to the relationship as the equivalence between SQP and Newton's method[67]. The (QP) subproblem is solved by utilizing the qpsolver[73] to call the open-source quadratic programming solver proxqp[74,75].

In this work, the $\ell_1$ merit function, defined as Eq.(23), is used to evaluate convergence criteria:

$$\phi_1(\mathbf{x}^k; \boldsymbol{\lambda}_{qp}^k, \boldsymbol{\mu}_{qp}^k) = f(\mathbf{x}^k) + (\boldsymbol{\rho}^k)^\top |\mathbf{h}(\mathbf{x}^k)| + (\mathbf{v}^k)^\top \mathbf{g}(\mathbf{x}^k)^-. \tag{23}$$

Where, $\mathbf{g}(\mathbf{x}^k)^- := \max(0, -\mathbf{g}(\mathbf{x}^k))$, and the penalty parameters are defined as follows[33]:

$$\boldsymbol{\rho}^k = \max\left(|\boldsymbol{\lambda}_{qp}^k|, \frac{\boldsymbol{\rho}^{k-1} + |\boldsymbol{\lambda}_{qp}^k|}{2}\right), \tag{24}$$

$$\mathbf{v}^k = \max\left(|\boldsymbol{\mu}_{qp}^k|, \frac{\mathbf{v}^{k-1} + |\boldsymbol{\mu}_{qp}^k|}{2}\right). \tag{25}$$

Given the $\ell_1$ merit function being not differentiable everywhere, the Eq.(26) describes the directional derivative of $\phi_1(\mathbf{x}^k; \boldsymbol{\lambda}_{qp}^k, \boldsymbol{\mu}_{qp}^k)$ along the direction $\mathbf{d}^k$ generated by the SQP subproblem.

$$D\left(\phi_1\left(\mathbf{x}^k;\boldsymbol{\lambda}_{qp}^k,\boldsymbol{\mu}_{qp}^k\right);\mathbf{d}^k\right) = \nabla f\left(\mathbf{x}^k\right)^\top \mathbf{d}^k - \left(\boldsymbol{\rho}^k\right)^\top \left|\mathbf{h}\left(\mathbf{x}^k\right)\right| - \left(\mathbf{v}^k\right)^\top \mathbf{g}\left(\mathbf{x}^k\right)^- \quad (26)$$

**Lemma 1.** *Let* $\mathbf{d}^k$, $\boldsymbol{\lambda}_{qp}^k$ *and* $\boldsymbol{\mu}_{qp}^k$ *be generated by problem* (QP). *Then the directional derivative of* $\phi_1\left(\mathbf{x}^k;\boldsymbol{\lambda}_{qp}^k,\boldsymbol{\mu}_{qp}^k\right)$ *in the direction* $\mathbf{d}^k$ *satisfies:*

$$D\left(\phi_1\left(\mathbf{x}^k;\boldsymbol{\lambda}_{qp}^k,\boldsymbol{\mu}_{qp}^k\right);\mathbf{d}^k\right) \leq 0 \quad (27)$$

**Proof.** By applying Taylor's theorem to $f$, $\mathbf{h}$ and $\mathbf{g}$, we obtain:

$$\begin{aligned}\phi_1\left(\mathbf{x}^k+\alpha\mathbf{d}^k;\boldsymbol{\lambda}_{qp}^k,\boldsymbol{\mu}_{qp}^k\right)-\phi_1\left(\mathbf{x}^k;\boldsymbol{\lambda}_{qp}^k,\boldsymbol{\mu}_{qp}^k\right) &= f\left(\mathbf{x}^k+\alpha\mathbf{d}^k\right)-f^k+\left|\mathbf{h}\left(\mathbf{x}^k+\alpha\mathbf{d}^k\right)\right|^\top\boldsymbol{\rho}^k-\left|\mathbf{h}^k\right|^\top\boldsymbol{\rho}^k+\left(\mathbf{v}^k\right)^\top\mathbf{g}\left(\mathbf{x}^k+\alpha\mathbf{d}^k\right)^- -\left(\mathbf{v}^k\right)^\top\mathbf{g}^{k-}\\ &\leq \alpha\nabla f^{k\top}\mathbf{d}^k+\gamma\alpha^2\left\|\mathbf{d}^k\right\|_2^2+\left(\left|\mathbf{h}^k+\alpha\left(\nabla\mathbf{h}^k\right)^\top\mathbf{d}^k\right|-\left|\mathbf{h}^k\right|\right)^\top\boldsymbol{\rho}^k+\left(\max\left(0,-\mathbf{g}^k-\alpha\left(\nabla\mathbf{g}^k\right)^\top\mathbf{d}^k\right)-\left(\mathbf{g}^{k-}\right)\right)^\top\mathbf{v}^k\end{aligned} \quad (28)$$

Where $f^k = f\left(\mathbf{x}^k\right)$, $\mathbf{h}^k = \mathbf{h}\left(\mathbf{x}^k\right)$, $\mathbf{g}^{k-} = \max\left(0,-\mathbf{g}\left(\mathbf{x}^k\right)\right)$, the positive constant $\gamma$ bounds the second-derivative terms in objective function and constraints, and $\alpha \in (0,1)$ is the step length.

According to Karush-Kuhn-Tucker (KKT) condition[67] of subproblem (QP), we have that $\mathbf{h}^k + \alpha\left(\nabla\mathbf{h}^k\right)^\top \mathbf{d}^k = 0$ and $\mathbf{g}^k + \alpha\left(\nabla\mathbf{g}^k\right)^\top \mathbf{d}^k \leq 0$, so for $\alpha \leq 1$ we have that

$$\begin{aligned}\phi_1\left(\mathbf{x}^k+\alpha\mathbf{d}^k;\boldsymbol{\lambda}_{qp}^k,\boldsymbol{\mu}_{qp}^k\right)-\phi_1\left(\mathbf{x}^k;\boldsymbol{\lambda}_{qp}^k,\boldsymbol{\mu}_{qp}^k\right) &\leq \alpha\nabla f^{k\top}\mathbf{d}^k+\gamma\alpha^2\left\|\mathbf{d}^k\right\|_2^2+\left((1-\alpha)\left|\mathbf{h}^k\right|-\left|\mathbf{h}^k\right|\right)^\top\boldsymbol{\rho}^k+\left((1-\alpha)\mathbf{g}^{k-}-\mathbf{g}^{k-}\right)^\top\mathbf{v}^k\\ &= \alpha\left[\nabla f^{k\top}\mathbf{d}^k-\left|\mathbf{h}^k\right|^\top\boldsymbol{\rho}^k-\left(\mathbf{g}^{k-}\right)^\top\mathbf{v}^k\right]+\gamma\alpha^2\left\|\mathbf{d}^k\right\|_2^2\end{aligned} \quad (29)$$

Similarly, the following lower bound can be obtained:

$$\phi_1\left(\mathbf{x}^k+\alpha\mathbf{d}^k;\boldsymbol{\lambda}_{qp}^k,\boldsymbol{\mu}_{qp}^k\right)-\phi_1\left(\mathbf{x}^k;\boldsymbol{\lambda}_{qp}^k,\boldsymbol{\mu}_{qp}^k\right) \geq \alpha\left[\nabla f^{k\top}\mathbf{d}^k-\left|\mathbf{h}^k\right|^\top\boldsymbol{\rho}^k-\left(\mathbf{g}^{k-}\right)^\top\mathbf{v}^k\right]-\gamma\alpha^2\left\|\mathbf{d}^k\right\|_2^2 \quad (30)$$

The directional derivative of $\phi_1$ in the direction $\mathbf{d}^k$ is given by Eq.(26) by taking the limit. Due to $\mathbf{d}^k$ satisfying the KKT condition of the subproblem (QP), we can reformulate the directional derivative as:

$$\begin{aligned}D\left(\phi_1\left(\mathbf{x}^k;\boldsymbol{\lambda}_{qp}^k,\boldsymbol{\mu}_{qp}^k\right);\mathbf{d}^k\right) &= -\mathbf{d}^{k\top}\mathbf{B}^k\mathbf{d}^k+\mathbf{d}^{k\top}\nabla\mathbf{h}^{k\top}\boldsymbol{\lambda}_{qp}^{k+1}+\mathbf{d}^{k\top}\nabla\mathbf{g}^{k\top}\boldsymbol{\mu}_{qp}^{k+1}-\left(\boldsymbol{\rho}^k\right)^\top\left|\mathbf{h}\left(\mathbf{x}^k\right)\right|-\left(\mathbf{v}^k\right)^\top\mathbf{g}\left(\mathbf{x}^k\right)^-\\ &\leq -\mathbf{d}^{k\top}\mathbf{B}^k\mathbf{d}^k+\left|\mathbf{h}^k\right|^\top\left(\boldsymbol{\lambda}_{qp}^{k+1}-\boldsymbol{\rho}^k\right)+\left(\mathbf{g}^{k-}\right)^\top\left(\boldsymbol{\mu}_{qp}^{k+1}-\mathbf{v}^k\right)\\ &\leq -\mathbf{d}^{k\top}\mathbf{B}^k\mathbf{d}^k+\left|\mathbf{h}^k\right|^\top\left(\boldsymbol{\lambda}_{qp}^{k+1}-\boldsymbol{\rho}^k\right)+\left(\mathbf{g}^{k-}\right)^\top\left(\boldsymbol{\mu}_{qp}^{k+1}-\mathbf{v}^k\right)\end{aligned} \quad (31)$$

Since $\mathbf{B}^k$ is positive definite, and according to Eq.(24) and Eq.(25), we have $\boldsymbol{\rho}^k \geq \boldsymbol{\lambda}_{qp}^{k+1}$ and $\mathbf{v}^k \geq \boldsymbol{\mu}_{qp}^{k+1}$, then the Lemma 1 has been proven. □

Lemma 1 indicates that if $\mathbf{d}^k$ is the KKT point of problem (QP) and $\mathbf{B}^k$ is positive definite, then the directional derivative decreases. This implies that in each iteration, either the infeasibility of the original problem is alleviated, or the value of the objective function is reduced, or both situations occur simultaneously.

The Armijo condition guarantees a sufficient decrease of the merit function:

$$\phi_1\left(\mathbf{x}^k + \alpha\mathbf{d}^k; \boldsymbol{\lambda}_{qp}^k, \boldsymbol{\mu}_{qp}^k\right) - \phi_1\left(\mathbf{x}^k; \boldsymbol{\lambda}_{qp}^k, \boldsymbol{\mu}_{qp}^k\right) < \alpha\eta D\left(\phi_1\left(\mathbf{x}^k; \boldsymbol{\lambda}_{qp}^k, \boldsymbol{\mu}_{qp}^k\right); \mathbf{d}^k\right). \tag{32}$$

Where $\eta \in (0, 0.5)$ is a hyper-parameter. However, it is possible that no step length $\alpha$ a can satisfy the Armijo condition, regardless of how small it is. Under such circumstances, we accept the last $\alpha$ when the maximum number of line searches has been reached.

For practical engineering problems, due to numerical stability and other factors, the KKT conditions as criteria may be too strict[76], making the algorithm difficult to achieve convergence. Therefore, based on previous research[33,77,78], we check whether the following criteria are met after every line search:

$$\left\|\mathbf{h}\left(\mathbf{x}^k + \alpha\mathbf{d}^k\right)\right\|_1 + \left\|\mathbf{g}\left(\mathbf{x}^k + \alpha\mathbf{d}^k\right)^-\right\|_1 < tol \tag{33}$$

$$\left|f\left(\mathbf{x}^k + \alpha\mathbf{d}^k\right) - f\left(\mathbf{x}^k\right)\right| < tol \tag{34}$$

$$\left\|\mathbf{d}^k\right\|_2 < tol \tag{35}$$

Where, $tol$ represents the convergence tolerance. The solution is considered optimal if Eqs.(33), (34) or Eqs.(33), (35) are satisfied. The satisfaction of Eq.(33) indicates that the optimal solution is feasible. The satisfaction of either Eq.(34) or Eq.(35) suggests that there is little potential for further decrease in the objective function.

## 4.2 Generalized framework

In this subsection, we will present a generalized framework that integrates SQP with the feasible path method to develop a derivative-based algorithm for differentiable machine learning models. The Generalized framework of SQP-driven feasible path algorithm for differentiable ML models is shown in the FIGURE 4.

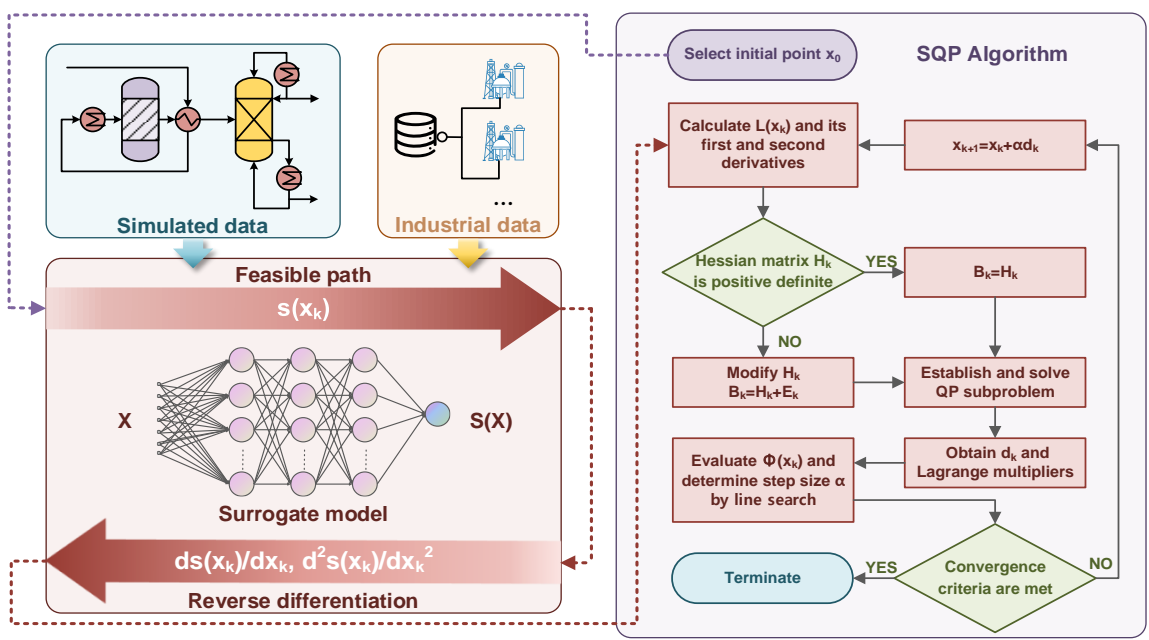

FIGURE 4 Generalized framework of SQP-driven feasible path algorithm for differentiable ML models.

For a surrogate-based optimization problem, the first step is to collect data from a rigorous simulator, a real-world factory or a combination of both types of data. The second step is to establish a surrogate model to approximate the relationship between input variables $\mathbf{x}$ and output $\mathbf{s}(\mathbf{x})$. The third step is to execute SQP algorithm. SQP-Driven feasible path algorithm for surrogate-based optimization problem as is shown in ALGORITHM 2.

ALGORITHM 2 SQP-Driven feasible path algorithm for surrogate-based optimization problem

| | |
|---|---|
| Step1: | Select a feasible initial point $\mathbf{x}^0$, initialize the $\boldsymbol{\lambda}_{qp}^0 = \boldsymbol{\rho}^0 = \boldsymbol{\mu}_{qp}^0 = \mathbf{v}^0 = 0$, $k=0$. Set hyper-parameter $\eta = 0.1$ and $tol = 10^6$. |
| Step2: | Calculate the output value of the surrogate model at $\mathbf{x}^k$, as well as the derivatives, $\mathbf{s}'$ and $\mathbf{s}''$, of the output $\mathbf{s}(\mathbf{x}^k)$ with respect to the input, which have been illustrated in Section 3. |
| Step3: | Calculate $\nabla f$, $\nabla \mathbf{h}$, $\nabla \mathbf{g}$ and $\mathbf{H}^k$ based on $\mathbf{s}(\mathbf{x}^k)$, $\mathbf{s}'$ and $\mathbf{s}''$. |
| Step4: | Check if $\mathbf{H}^k$ is positive definite. If so, set $\mathbf{B}^k = \mathbf{H}^k$; otherwise, set $\mathbf{B}^k = \mathbf{H}^k + \mathbf{E}^k$ as is shown in Section4. |
| Step5: | Establish and solve (QP) subproblem, then obtain the direction $\mathbf{d}^k$ and Lagrangian multiplier $\boldsymbol{\lambda}_{qp}^k$ and $\boldsymbol{\mu}_{qp}^k$. Update $\boldsymbol{\rho}^k$ and $\mathbf{v}^k$. |
| Step6: | Calculate merit function $\phi$ and directional derivative $D$. |
| Step7: | Use the backtracking method to search the step size $\alpha$ satisfying the Armijo condition, with an iteration multiple of 0.618 and a maximum iteration count of 10. |
| Step8: | Check Eqs.(33), (34) or Eqs.(33), (35) are satisfied. If so, terminate the algorithm; otherwise, $\mathbf{x}^{k+1} = \mathbf{x}^k + \alpha \mathbf{d}^k$, $k = k+1$ and back to Step2. |

The original problem is decomposed into inner and outer loops by feasible path method. To reduce the problem size, the forward mode is employed to calculate the value of the iteration point, instead of solving a group of algebraic equations generated by ML models. The backward mode is used to calculate the derivatives to formulate QP subproblems and indicate the direction for the next iteration. The proposed framework is a universal, efficient, and effective algorithm for optimization problems embedded with differentiable machine learning models.

## 5. Case study

Given the challenges in acquiring data from real-world production processes, we initially developed and optimized surrogate models for six test functions with known minimum values. This part served as a validation study to ensure the effectiveness of our framework. Following this validation, we applied our methodologies to two case studies based on Aspen Plus simulations, where we also established and optimized surrogate models. Throughout these case studies, we consistently used MLPs to implement all surrogate models with Swish activation function as is shown in Eq.(36). Latin Hypercube Sampling (LHS) is employed to construct data sets, which ensures a robust and representative distribution of samples across the input space. This approach not only guarantees consistency but also facilitates a more accurate comparison between different scenarios.

$$Swish(x) = \frac{x}{1+e^{-x}} \tag{36}$$

### 5.1 Test functions

TABLE 1 Test functions

| Function name | Dimens-ion | Expression | Minimal value | Minimum point |
|---|---|---|---|---|
| Sphere | 10 | $\sum_{i=1}^{10} x_i^2$ | 0 | (0,0,…,0) |
| Quadratic | 10 | $\sum_{i=1}^{10} x_i^2 + \sum_{i=1}^{9}(x_i - x_{i+1})^2$ | 0 | (0,0,…,0) |
| Six-hump camel | 2 | $\left(4 - 2.1x_1^2 + \frac{x_1^4}{3}\right)x_1^2 + x_1 x_2 + 4(-1 + x_2^2)x_2^2$ | -1.03 | (0.09,-0.71) (-0.09,0.71) |

| Schaffer NO.2 | 2 | $0.5 + \dfrac{\sin^2\left(x_1^2 - x_2^2\right) - 0.5}{\left[1 + 0.001\left(x_1^2 + x_2^2\right)\right]^2}$ | 0 | (0,0) |
| Griewank | 5 | $1 + \dfrac{1}{4000}\sum_{i=1}^{5} x_i^2 - \prod_{i=1}^{5} \cos\left(\dfrac{x_i}{\sqrt{i}}\right)$ | 0 | (0,0,0,0,0) |
| Ackley | 5 | $-20\exp\left(-0.2\sqrt{\dfrac{1}{5}\sum_{i=1}^{5} x_i^2}\right) - \exp\left(\dfrac{1}{5}\sum_{i=1}^{5}\cos(2\pi x_i)\right) + 20 + e$ | 0 | (0,0,0,0,0) |

Five typical functions as shown in TABLE 1 are selected as test functions. Except for the Sphere function and a special form of Quadratic function, which are high-dimensional convex functions, the other four functions are nonconvex and possess multiple local minimum points. The sampling range of the six test functions is selected as [-2,2], with a training set sample size of 2000, a validation set of 50, and a test set of 200. The performance of the six MLPs, as illustrated in FIGURE 5, demonstrates their capability to serve as effective and accurate surrogates for the original test functions.

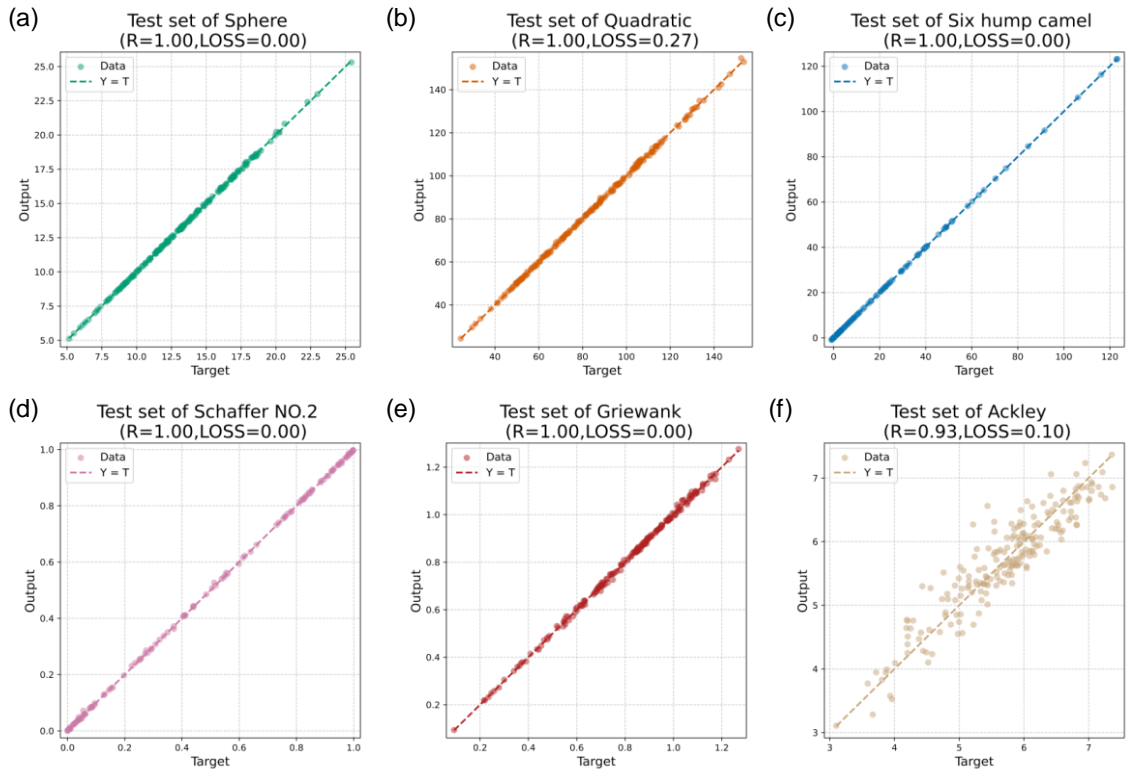

FIGURE 5 The parity plots over test set with pearson correlation coefficient (R) and testing loss (LOSS) to evaluate the accuracy of surrogate models for (a) Sphere function, (b) Quadratic function, (c) Six hump camel function, (d) Schaffer function NO.2, (e) Griewank function and (f) Ackley function.

FIGURE 6 shows the convergence curves of the SQP-driven feasible path algorithm for solving six surrogate models with all initial points at $x_i = 2$. $\|\hat{x}^* - x^*\|_2$ represent the distance between solved optimal solution $\hat{x}^*$ and true global optimal solution $x^*$, which in all six cases does not exceed 0.1. According to the true function value at the obtained solution and the global minimum value, the results indicate that the surrogate model combined with the SQP-driven feasible path algorithm successfully converges to the optimum in an extremely short time.

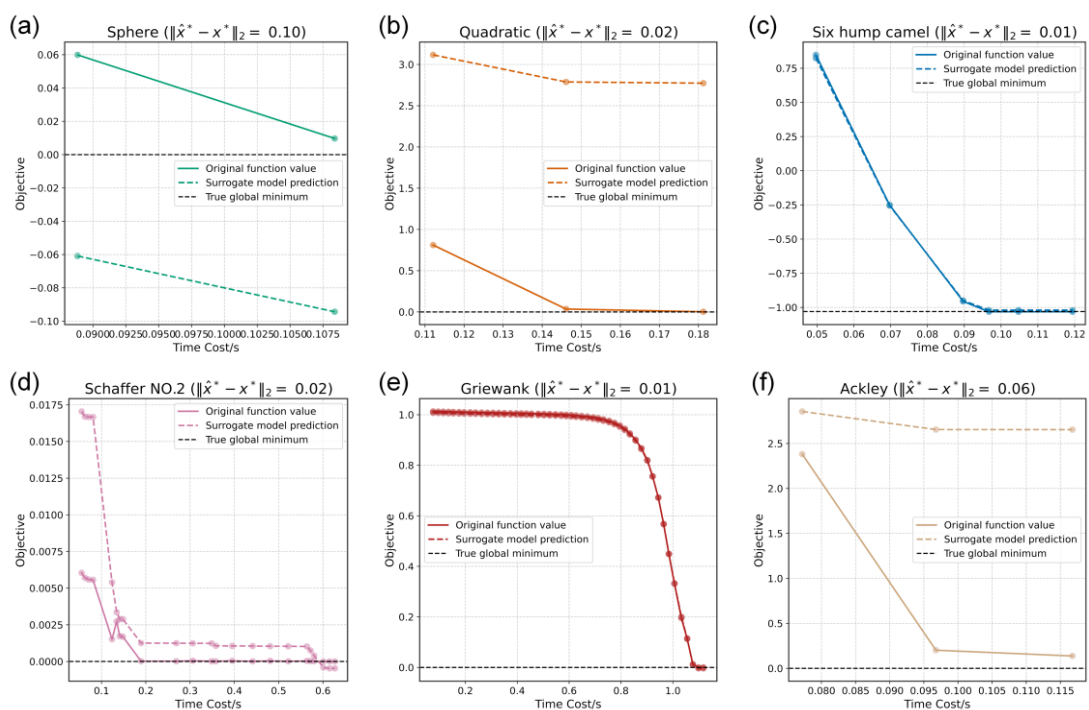

FIGURE 6 The convergence curve of surrogate models for (a) Sphere function, (b) Quadratic function, (c) Six hump camel function, (d) Schaffer function NO.2, (e) Griewank function and (f) Ackley function with distance between solved optimum and true global optimum.

## 5.2  Extractive distillation case

In this subsection, we tested an example of extractive distillation to separate toluene form n-heptane using solvent phenol provided by Ma et al[54], as is shown in FIGURE 7. The process mainly consists of two columns. The first column C1 is used for separate n-heptane, and the second column C2 is used for separate toluene form solvent phenol.

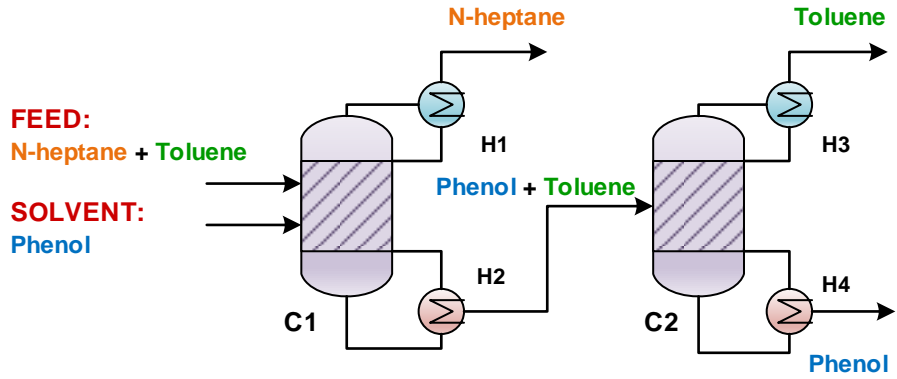

FIGURE 7 The process flow diagram for extractive distillation to separate toluene form n-heptane using solvent phenol in Ma et al.'s work[54].

The optimization model mainly achieves the minimum operating cost by adjusting the amount of extractant, reflux ratio, and distillate flow rate, as is shown in Eq.(37). The purity requirement for both products is not less than 0.97.

$$\begin{aligned}
&\min_{F,r_1,D_1,r_2,D_2 \in \mathbb{R}} && c_{hu}(Q_2+Q_4) - c_{cu}(Q_1+Q_3) \\
&\text{s.t.} && s(F,r_1,D_1,r_2,D_2) \rightarrow (x_n,x_t,Q_1,Q_2,Q_3,Q_4) \\
& && 0.97 \leq x_n \leq 1 \\
& && 0.97 \leq x_t \leq 1
\end{aligned} \quad (37)$$

Where, $c_{hu}$ and $c_{cu}$ is the energy cost parameters for cold and hot utilities, set to 20 and 80 \$·kW$^{-1}$·yr$^{-1}$ [79]; $Q_1$ and $Q_3$ are the condenser duty of C1 and C2, $Q_2$ and $Q_4$ are the reboiler duty of C1 and C2, kW; $F$ is the mass flow of solvent phenol, kg·hr$^{-1}$; $r_1$ and $r_2$ are the reflux ratio of C1 and C2; $D_1$ and $D_2$ are the distillate flow of C1 and C2; $x_n$ and $x_t$ are the mole fraction of n-heptane in C1 distillate and toluene in C2 distillate; $s(\cdot)$ is the surrogate model for extractive distillation process.

According to Ma et al.[54], surrogate models constructed using either ReLU neural networks or ALAMO (Automated Learning of Algebraic Models for Optimization) exhibit relative errors ranging from 0.2% to 25%. Despite employing complex basis functions, these models still fail to fully capture all behaviors of the process system due to inherent limitations in their nonlinear representation capabilities.

However, in this work, by introducing the feasible path method, which eliminates a significant number of intermediate variables, we achieve end-to-end inference and differentiation of the surrogate model. Consequently, we could pay more attention on whether the model can accurately capture the behavior of the system, without worrying about the complexity of the network structure.

The sampling numbers for the training set, validation set, and test set were set to 2000, 200, and 1000, with 1014, 95, and 511 successful points, respectively. The surrogate model in this case also constructed with MLP, who possesses four hidden layers with hidden layers 1-2 having 100 neurons and hidden layers 3-4 having 120 neurons. FIGURE 8 shows the performance of constructed model on test set. The relationship between the predicted values and the target values for all outputs depicted in FIGURE 8 closely approximates the perfect linear correspondence, since each graph exhibits a correlation coefficient of 1.00 and relative errors are no greater than 0.36%. This indicates excellent predictive performance of the surrogate model. Such results are highly encouraging for model validation and instill strong confidence in the models' reliability.

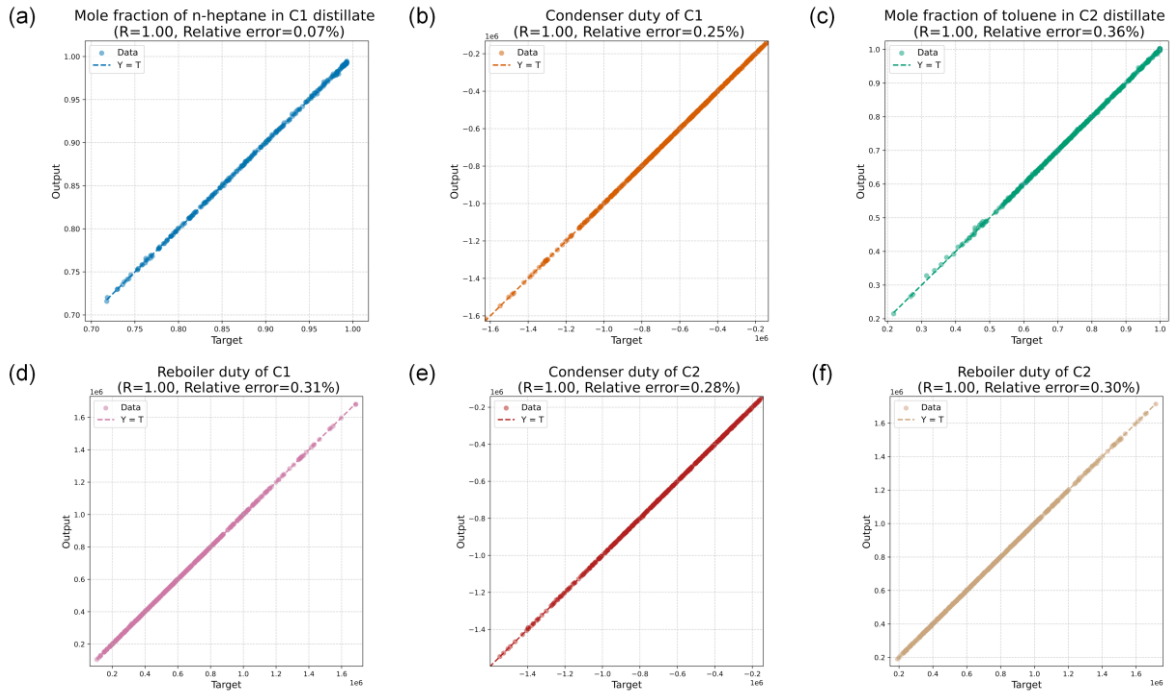

FIGURE 8 The parity plots over test set with pearson correlation coefficient (R) and relative error on test set to evaluate the accuracy of surrogate model's outputs, specifically (a) mole fraction of n-heptane in C1 distillate, (b) mole fraction of toluene in C2 distillate, (c) condenser duty of C1, (d) reboiler duty of C1, (e) condenser duty of C2, and (f) reboiler duty of C2 of extractive distillation process.

The accurate surrogate models are the cornerstone of high-quality optimization solution. Although Ma et al.[54] globally optimized the constructed surrogate ReLU NN or algebraic models, the unreliable model ultimately led to infeasibility of purity constraints. Therefore, in this work, we must not only adhere to the constraints of the surrogate optimization model but also verify the constraints derived from the original Aspen Plus simulation model. The convergence curve is presented in FIGURE 9. This figure demonstrates that, apart from the failure points in the Aspen Plus simulation, there is good consistency between the surrogate MLP model and the Aspen simulation results, with all purity constraints being satisfied. Consequently, we can confidently assert that the optimization has been successful.

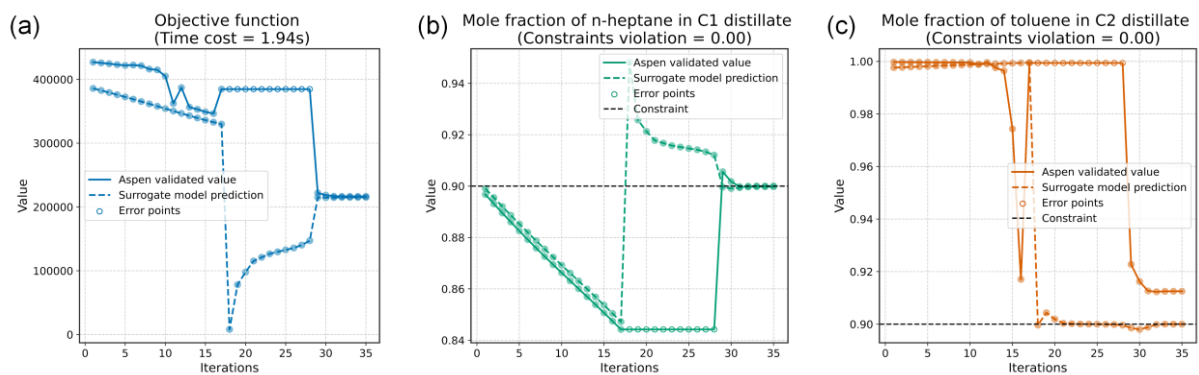

FIGURE 9 The convergence curve of (a) objective using surrogate model pridiction and Aspen validated value, and surrogate model pridiction as well as Aspen validated value for the mole fraction of (b) n-heptane in C1 distillate and (c) toluene in C2 distillate of extractive distillation process.

## 5.3  Process of absorbing $CO_2$ from biogas

Removing $CO_2$ from biogas is an important step in producing biofuels as it can increase the concentration of $CH_4$, thereby enhancing the quality of biogas as an energy carrier. Purified biogas, commonly known as biomethane or renewable natural gas (RNG), can be used for various purposes, such as direct combustion to generate heat or electricity, as vehicle fuel, or injection into natural gas pipelines.

In this case, first introduced by Xu et al.[80], we remove $CO_2$ with n-methyldiethanolamine (MDEA). The reactions between $CO_2$ and MDEA takes place in the absorber, as illustrated in Eq. (38). Then, the carbonic acid undergoes thermal decomposition in the desorber, as illustrated in Eq.(39).

$$CO_2 + R_2NCH_3 + H_2O \longrightarrow R_2NCH_3H^+ + HCO_3^- \tag{38}$$

$$R_2NCH_3H^+ + HCO_3^- \xrightarrow{\Delta} R_2NCH_3 + CO_2 \uparrow + H_2O \qquad (39)$$

The desorbed MDEA is circulated in the system as the regenerated solvent. The products CH$_4$ and CO$_2$ will carry a small amount of water vapor and MDEA, so it is necessary to supplement the solvent appropriately. Therefore, we employ Python to externally call Aspen's Application Programming Interface (API) for direct iteration, ensuring that the loss remains consistent with the replenished solvent.

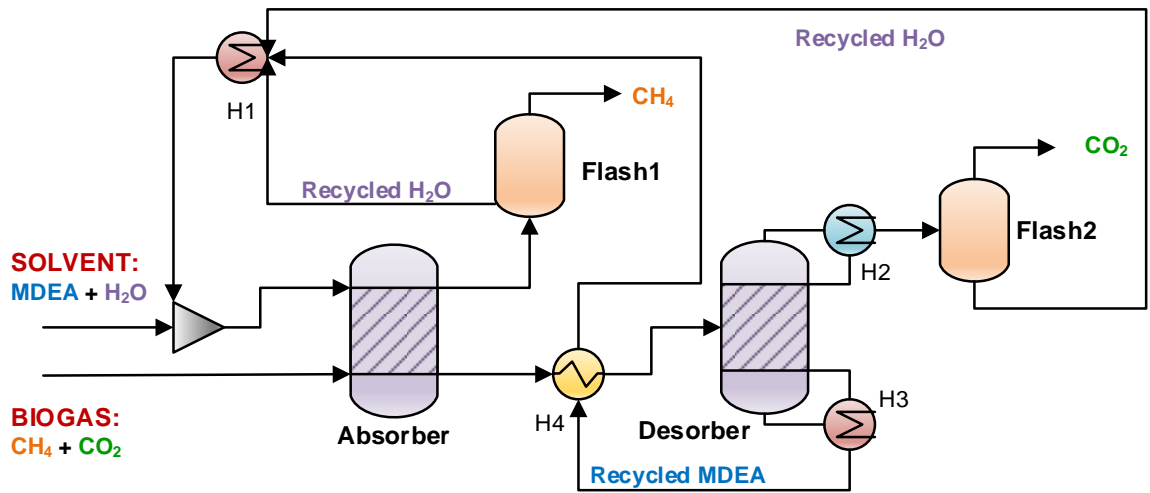

FIGURE 10 The process flow diagram for absorbing CO$_2$ in biogas using aqueous MDEA solution.

The optimization model is presented in Eq.(40). The objective is to minimize annual total cost, including utility costs, equipment investment costs, and supplementary MDEA costs. The process stipulates that the purified biogas must achieve a minimum CH$_4$ purity of 0.97, while the CO$_2$ content in the regenerated MDEA solution should be less than 0.1.

$$\min_{\substack{F^R_{MDEA} \in \mathbb{R} \\ N_A, N_D \in \mathbb{Z}}} \quad HrF^S_{MDEA} - c_{cu}Q_c + c_{hu}Q_h + \left[ c_{MDEA}F^R_{MDEA} + c_s(N_A + N_D) \right]/Yr$$

$$\begin{aligned} s.t. \quad & s\left(F^R_{MDEA}, N_A, N_D\right) \rightarrow \left(x_{CH_4}, x_{CO_2}, F^S_{MDEA}, Q_c, Q_h\right) \\ & 0.97 \leq x_{CH_4} \leq 1 \\ & 0 \leq x_{CO_2} \leq 0.1 \\ & 2.7 \leq F^R_{MDEA} \leq 5 \\ & 3 \leq N_A \leq 5 \\ & 3 \leq N_D \leq 11 \end{aligned} \quad (40)$$

Where, $Hr$ is the annual operating time, set to 8000 hours; $Yr$ is the period of depreciation, set to 10 years; $c_{hu}$ and $c_{cu}$ is the energy cost parameters for cold and hot utilities, set to 20 and 80 $·kW$^{-1}$·yr$^{-1}$ [79]; $c_{MDEA}$ is the price of MDEA, set to 13000 $·ton$^{-1}$; $c_s$ is the price for each column stage, set to 8508 $ [81,82]; $Q_c$ and $Q_h$ are cooling and heating utilities of the whole system, kW; $F^R_{MDEA}$ is the mass flow of recycled MDEA, ton·hr$^{-1}$; $F^S_{MDEA}$ is the mass flow of supplementary MDEA, ton·hr$^{-1}$; $N_A$ and $N_D$ are the number of stages of absorber and desorber; $x_{CH_4}$ and $x_{CO_2}$ are the mass fraction of $CH_4$ and $CO_2$ in system outlets; $s(\cdot)$ is the surrogate model for absorption process. Since the proposed algorithm is intended to solve NLP problems, we apply it to solve the relaxed subproblem of Eq.(40), where $N_A, N_D \in \mathbb{R}$.

The surrogate was model also trained as MLP, who possesses six hidden layers with 1-3 layers having 100 neural and 4-6 layer having 200 neural. The surrogate model in this case also constructed with MLP, who possesses four hidden layers with hidden layers 1-2 having 100 neurons and hidden layers 3-4 having 120 neurons. The sampling numbers for the training set, validation set, and test set were set to 2000, 100, and 200. The performance on test set as is shown in FIGURE 11, which show the great consistency between surrogate model and first principles model.

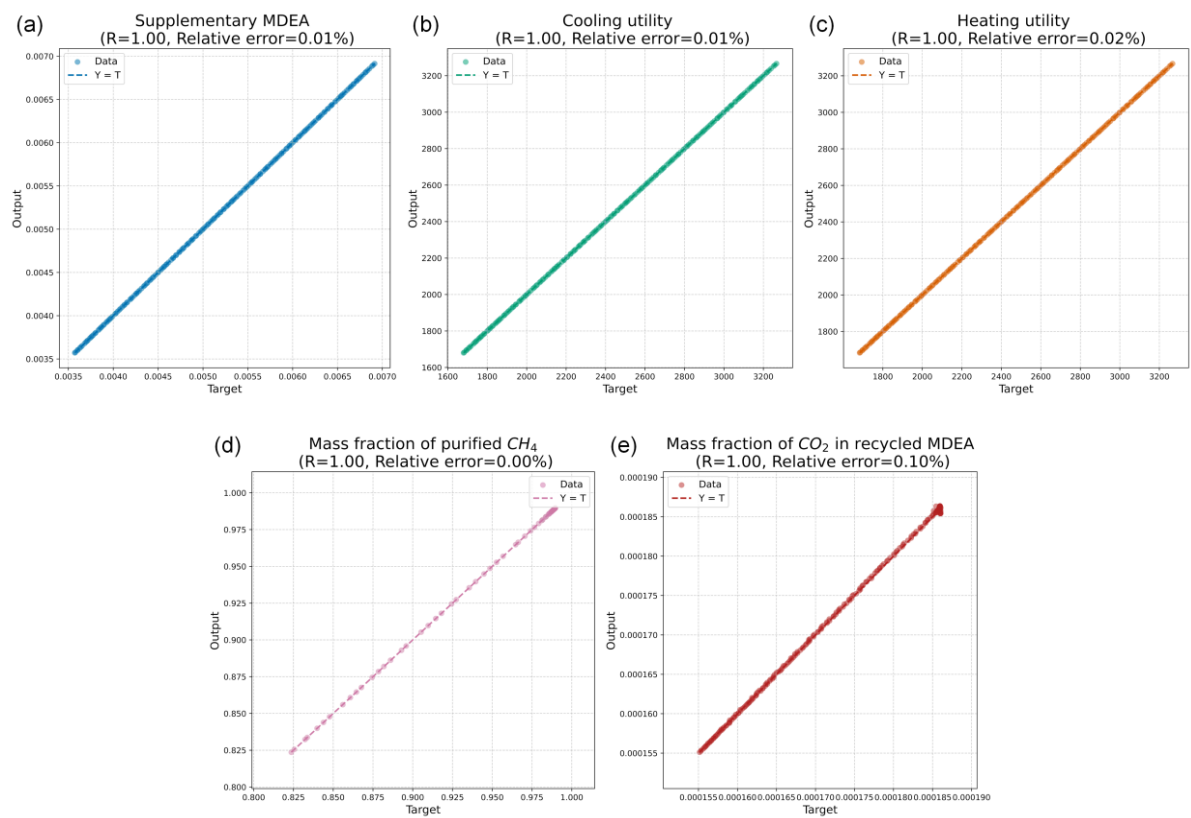

FIGURE 11 The parity plots over test set with pearson correlation coefficient (R) and relative error on test set to evaluate the accuracy of surrogate model's outputs, specifically (a) supplementary MDEA, (b) cooling utility, (c) heating utility, (d) mass fraction of purified $CH_4$, (e) mass fraction of $CO_2$ in recycled MEDA, and (f) reboiler duty of C2 of process for absorbing $CO_2$ from biogas.

The convergence curve of the relaxed NLP problem as is shown in FIGURE 12, where the curve for surrogate model prediction and Aspen validated value are perfectly coinciding. Usually, the iteration points of relaxation problems are not feasible for primal problem. While fortunately in this case, we reached the optimal solution in two iterations, and each point is a feasible solution to the primal problem. This case highlights the promising potential of the proposed algorithm in surrogate-based MINLP problems, such as process synthesis, reaction-distillation-heat-exchange-network coupled optimization, etc. The algorithm can effectively serve as a submodule within the MINLP framework, specifically for solving NLP subproblems.

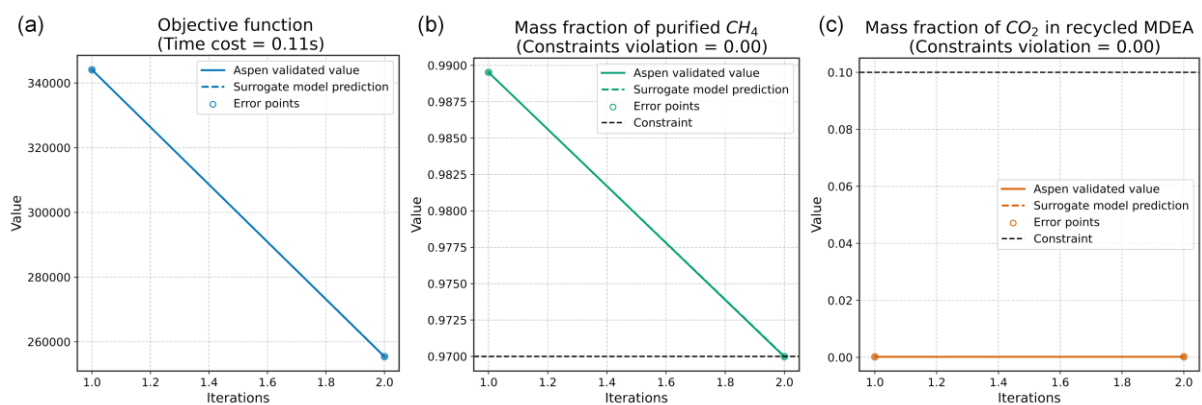

FIGURE 12 The convergence curve of (a) objective using surrogate model pridiction and Aspen validated value, and surrogate model pridiction as well as Aspen validated value for the mole fraction of (b) mass fraction of purified CH$_4$, and (c) mass fraction of CO$_2$ in recycled MEDA of process for absorbing CO$_2$ from biogas.

## 6. Conclusion

A generalized framework driven by feasible path SQP algorithm is proposed for solving simulation-based problems with differentiable ML models serving as surrogates. By incorporating the concept of a feasible path, the framework manages to avoid the proliferation of intermediate variables that typically arise from the algebraic formulations of machine learning models, while successfully preserving the rapid inference capabilities of ML models. The first- and second-order derivatives of the ML model outputs with respect to the inputs are derived, enabling their integration into any derivative-driven algorithms. The results of the case studies demonstrate the efficiency and effectiveness of the proposed algorithm, as all six test functions were successfully globally optimized, and the optimization step was completed within 2 seconds for each case. The Aspen simulation verification confirms the consistency between the surrogate models and the rigorous models, ensuring that the solutions adhere to the constraints of the underlying systems.

Building on this foundation, the proposed framework opens up promising avenues for future research in both NLP and MINLP problems. It can be integrated into global optimization algorithms for NLP to address surrogate-based upper or lower bounds or serve as a submodule in MINLP algorithms for solving NLP subproblems. Moreover, the robustness and accuracy of the surrogate models should be strengthened in future work, since different surrogates can lead to different optima, and only reliable surrogates can ensure successful optimization. Additionally, models with high differentiation complexity but common usage, such as deep neural networks and recurrent neural networks, remain challenging to incorporate into this framework. We will continue to explore these challenges in our future research.

**Supporting information**

Additional supporting information can be found online in the Supporting Information section at the end of this article. The raw data set sampled from 3 cases are provided in Data A as a .zip file. The Aspen Plus simulations of case 2 and 3 are provided in Data B as a .zip file. The neural network surrogate models for 3 cases are provided in Data C as a .zip file. An available example code of proposed algorithm is provided at https://github.com/ZixuanChang/MLSQP .